\documentclass[11pt]{article}
\usepackage{epsf,amssymb,latexsym,amsmath}
\usepackage{graphicx}
\usepackage{tikz}

\newcommand{\proof}{\noindent{\bf Proof.\ }}
\newcommand{\qed}{\hfill $\square$ \bigskip}

\newtheorem{theorem}{\bf Theorem}[section]
\newtheorem{corollary}[theorem]{\bf Corollary}
\newtheorem{lemma}[theorem]{\bf Lemma}

\newtheorem{definition}[theorem]{\bf Definition}

\begin{document}

\title{The domination number of plane triangulations}

\author{
Simon \v Spacapan\footnote{ University of Maribor, FME, Smetanova 17,
2000 Maribor, Slovenia. e-mail: simon.spacapan @um.si.
}}
\date{\today}

\maketitle

\begin{abstract}
\noindent We introduce a class of plane graphs called weak near-triangulations, and 
 prove that this class  is closed under certain graph operations. 
Then we use the properties of weak near-triangulations to prove that every plane triangulation on $n>6$ vertices has a dominating set of size at most $17n/53$. 
This improves the 
bound $n/3$ obtained by Matheson and Tarjan. 
\end{abstract}



\section{Introduction}
 A dominating set in a graph $G=(V,E)$ is a set $D\subseteq V$, such that every vertex in $V\setminus D$ is adjacent to a vertex in $D$. The domination number of 
$G$, denoted as $\gamma(G)$, is the minimum size of a dominating set in $G$. 
A plane triangulation is a connected simple plane graph $G$, such that every face of $G$ is triangular. A plane graph $G$ is a  near-triangulation, if $G$ is 2-connected and  every  face of $G$ is triangular, except possibly the unbounded face.


Matheson and Tarjan conjectured in  \cite{mathesontarjan} 
that every sufficiently large  plane triangulation  has a dominating set of cardinality at most $n/4$, where $n$ is the number of vertices of the triangulation. 
They proved that the domination number of every near-triangulation, and therefore also every triangulation,  is at most $n/3$. 

\begin{theorem}\cite{mathesontarjan}\label{mathesontarjan}
Every near-triangulation has a 
3-coloring such that each color class is a dominating set. 
\end{theorem}
We mention that Theorem \ref{mathesontarjan} provides the only known upper bound for $\gamma(G)/n$
in the class of plane triangulations. Several results for sublasses of plane triangulations are known. For example, the conjecture of Matheson and Tarjan was confirmed in  \cite{erika}  for plane triangulations with maximum degree 6. In \cite{liu} the authors  improve the bound $n/4$ obtained in \cite{erika}, by proving that there exists a constant $c$ such that  $\gamma(G)\leq n/6+c$ 
 for every plane triangulation $G$ with maximum degree 6. 
Another result is given in \cite{plummer1} where the authors prove that $\gamma(G)\leq \max\{\lceil\frac {2n}{7}\rceil,\lfloor\frac{5n}{16}\rfloor\}$ 
for every 4-connected plane triangulation $G$. 
The domination number of outerplanar triangulations is considered in \cite{campos} and  \cite{tokunaga}, where the authors  independently prove  that an outerplaner triangulation 
with $n$ vertices and $t$ vertices of degree 2 has a dominating set of cardinality at most $(n+t)/4$. This result is further improved in \cite{li}. 

Theorem \ref{mathesontarjan} was extended to tringulations on the projective plane, the torus and the Klein bottle in \cite{plummer} and in \cite{honjo}. It is proved that every triangulation 
on any of  these surfaces  has a dominating set of cardinality  at most $n/3$. These results are further generalized in \cite{furuya}, where it is proved that every triangulation on a closed surface has a dominating set of cardinality at most $n/3$. 

We also mention that the domination number of planar graphs with small diameter is studied in \cite{mac}, \cite{henning1} and \cite{henning2}. It is proved that every sufficiently large planar 
 graph with diameter 3 has domination number at most 6.

We approach the problem of finding the smallest constant $c$, such that $\gamma(G)\leq cn$, for every sufficiently large plane triangulation. In the following section we introduce the class of weak near-triangulations. We define  reducibility of a weak near-triangulation, and prove that every  weak near-triangulation $G$ is reducible, except if all blocks of $G$ are outerplaner or 
have exactly 6 vertices. 
This result is then used to prove that  every plane triangulation on $n>6$ vertices has a dominating set of cardinality at most 17n/53.


\section{The domination number of triangulations}
We refer the reader to \cite{diestel} and \cite{mt} for a complete overwiev of definitions and terminology that we use. 
In this article we consider only simple graphs with no multiple edges or loops. Let $G$ be a plane graph.  
An edge $e$   of $G$ is {\em incident} to a face $F$ of $G$ if $e$  is contained in the boundary of $F$. Similarly we define incidence of a verex $x$  and a face $F$.  
 Vertices and edges  incident to the unbounded face of $G$ are called  {\em external vertices} and {\em  external edges}, respectively.  
If a vertex (an edge) is not an external vertex (edge), then 
 it is called an  {\em internal vertex} ({\em  internal edge}). 
 A {\em triangle} is a cycle on three vertices. If $a,b$ and $c$ are vertices of a triangle, we denote this triangle by $abc$, and we say that $a,b$ and $c$ are {\em contained} in the triangle $abc$.  A {\em facial triangle} of a plane graph $G$ is a triangle whose interior is a face of $G$.  A face bounded by a triangle is called a {\em triangular face}. A connected plane graph $G$  is a {\em triangulation}  if all faces of $G$ are triangular, and $G$ is a  {\em near-triangulation} if it is 2-connected and  every bounded face of $G$ is triangular. 
A {\em block} of a graph $G$ is a maximal connected subgraph of $G$ without a cutvertex.



\begin{definition}
A plane graph $G$ is a weak near-triangulation (WNT) if every bounded face of $G$ is triangular  and every vertex of $G$ is contained in a triangle.
\end{definition}

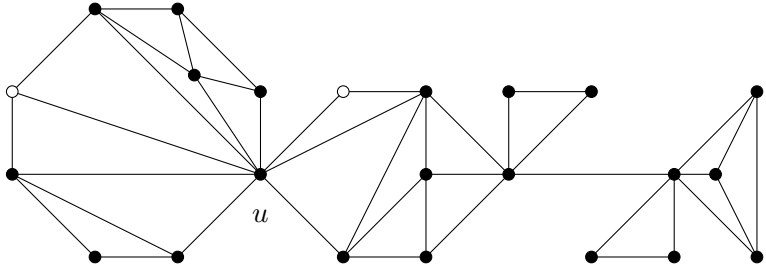
\begin{figure}[htb] 
\begin{center}
\begin{tikzpicture}[scale=1.1]

\filldraw (0,0) circle (2pt);
\filldraw (0,1) circle (2pt);
\filldraw (-1,2) circle (2pt);
\filldraw (-2,2) circle (2pt);

\filldraw (-3,0) circle (2pt);
\filldraw (-2,-1) circle (2pt);
\filldraw (-1,-1) circle (2pt);
\filldraw (-0.8,1.2) circle (2pt);

\filldraw (2,1) circle (2pt);
\filldraw (2,0) circle (2pt);
\filldraw (3,0) circle (2pt);

\filldraw (1,-1) circle (2pt);
\filldraw (2,-1) circle (2pt);

\filldraw (3,1) circle (2pt);
\filldraw (4,1) circle (2pt);
\filldraw (5,0) circle (2pt);
\filldraw (6,1) circle (2pt);

\filldraw (6,-1) circle (2pt);
\filldraw (5.5,0) circle (2pt);

\filldraw (4,-1) circle (2pt);
\filldraw (5,-1) circle (2pt);



\draw  (0,0) to (-1,-1);
\draw  (-1,-1) to (-2,-1);
\draw  (-2,-1) to (-3,0);
\draw  (-3,0) to (-3,1);
\draw  (-2,2) to (-3,1);
\draw  (-2,2) to (-1,2);
\draw  (-1,2) to (0,1);
\draw  (0,1) to (0,0);
\draw  (0,0) to (-3,0);
\draw  (0,0) to (-3,1);
\draw  (0,0) to (-2,2);
\draw  (-3,0) to (-1,-1);
\draw  (-0.8,1.2) to (0,0);
\draw  (-0.8,1.2) to (0,1);
\draw  (-0.8,1.2) to (-1,2);
\draw  (-0.8,1.2) to (-2,2);

\draw  (0,0) to (2,1);
\draw  (2,1) to (1,-1);
\draw  (0,0) to (1,1);
\draw  (1,1) to (2,1);
\draw  (2,1) to (3,0);
\draw  (3,0) to (2,-1);

\draw  (0,0) to (1,-1);
\draw  (1,-1) to (2,-1);

\draw  (2,0) to (2,1);
\draw  (2,0) to (3,0);
\draw  (2,0) to (2,-1);
\draw  (2,0) to (1,-1);

\draw  (3,0) to (5,0);
\draw  (3,0) to (3,1);
\draw  (3,1) to (4,1);
\draw  (3,0) to (4,1);

\draw  (5,0) to (5.5,0);
\draw  (5,0) to (6,1);
\draw  (5,0) to (6,-1);
\draw  (6,1) to (6,-1);
\draw  (6,1) to (5.5,0);
\draw  (6,-1) to (5.5,0);

\draw  (4,-1) to (5,-1);
\draw  (4,-1) to (5,0);
\draw  (5,0) to (5,-1);

\path node at (0,-0.5) {$u$};

\filldraw[fill=white] (-3,1) circle (2pt);
\filldraw[fill=white] (1,1) circle (2pt);

\end{tikzpicture}
\caption{A weak near-triangulation $G$. \label{sibkaskorajtriangulacija}} 
\end{center}
\end{figure}

We note that an empty graph is a weak near-triangulation, and that in the above definition there are no assumptions about the connectivity of $G$, hence $G$ may be disconnected. Note also that  the definition is equivalent to  
the following: $G$ is a weak near-triangulation if every bounded face of $G$ is triangular and every vertex of $G$ is contained in a facial triangle. 

In \cite{diestel}  (see Lemma 4.2.2., p. 91, and Lemma 4.2.6, p. 93) the following two results are given. 
\begin{lemma}\label{prva}
Let $G$ be a plane graph and $e$ an edge of $G$. If $e$ lies on a cycle $C\subseteq G$, then $e$ is incident to exactly one face $F$ of $G$, such that $F$ is contained in the interior of $C$. 
\end{lemma}

\begin{lemma}\label{druga}
In a 2-connected plane graph, every face is bounded by a cycle. 
\end{lemma}

We use Lemma \ref{prva} and \ref{druga} to prove the following. 

\begin{lemma} \label{blok}
Every block of a weak near-triangulation is either a near-triangulation or a $K_2$.
\end{lemma}
\proof 
Let $B$ be a block of a weak near-triangulation $G$. Since every vertex of $G$ is contained in a triangle of $G$, we find that $B$ has more than one vertex. 
Suppose that $B$ is a block on more than two vertices. Then $B$ is 2-connected, and so by Lemma \ref{druga} every face of $B$ is bounded by a cycle. 
Let $F$ be a bounded face of $B$ and $C$ the cycle bounding $F$ ($F$  lies in the interior of $C$). 
To prove that $B$ is a near-triangulation, we have to prove that $F$ is triangular.  

Let $e=uv$ be any edge of $C$. Since $C$ is a cycle of $G$ we find, by Lemma \ref{prva}, 
that there is a face $F'$ of $G$ contained in the interior of $C$ incident to $e$. 
Since $F'$ is a bounded face of $G$, it has to be triangular (by the definition of a weak near-triangulation). If $F'=F$ we are done. Suppose that $F'\neq F$, and suppose that $F'$ is bounded by the triangle $uvx$. Then $x$ is not in $B$, for otherwise $F$ and $F'$ are faces of $B$ contained in the interior of $C$, both incident to $e$ (contradictory to Lemma \ref{prva}). So $x\notin B$ and therefore $B$ is not a maximal connected subgraph of $G$ without a cutvertex (we may add $x$ to $B$), a contradiction.
\qed

Observe also that every endblock (leaf block) of a weak near-triangulation is a near-triangulation. 
\begin{lemma}
Every weak near-triangulation has a 
3-coloring such that each color class is a dominating set. 
\end{lemma}

\proof 
Let $G$ be a weak near-triangulation. Delete all bridges of $G$ and call the obtained graph $G'$. 
Every block of $G'$ is a near-triangulation, and two distinct blocks have at most one common vertex. 
To obtain the desired 3-coloring  of a connected component of $G'$ we use induction on number of blocks. 
In induction step we delete an endblock (except the cutvertex) and color the obtained graph with 3 colors according to the induction hypothesis.   
Then use Theorem \ref{mathesontarjan} to obtain a coloring of the deleted endblock, if needed permute the color classes in the endblock (so that the cutvertex gets 
the color that it already has).  \qed

A straightforward corollary is the following. 

\begin{corollary}\label{n3weak}
Every weak near-triangulation on $n$ vertices has domination number at most $n/3$. 
\end{corollary}

For a plane graph $G$ and $X\subseteq V(G)$ we denote by $G-X$ the graph obtained from $G$ by deleting the vertices in $X$ and all edges incident to a vertex in $X$. 
If $X=\{u\}$ is a singleton, we write $G-u$ instead of $G-\{u\}$. If $H$ is a subgraph of a 
plane graph $G$, then we assume that the embedding of $H$ in the plane is that given by $G$. This in particular applies for $H=G-X$.    Let $N[x]$ denote  
 the set of vertices that are either adjacent or equal to $x$, $N[x]$ is called  the {\em closed 
neighborhood} of $x$. 

\begin{definition}
Let $G$ be a weak near-triangulation.  We say that $G$ is reducible if there exists a set $D\subseteq V(G)$ and a vertex $x\in D$  with the following properties:
\begin{itemize}
\item[(i)] $D\subseteq N[x]$
\item[(ii)] $|D|\geq 4$
\item[(iii)]  $G-D$ is a weak near-triangulation.
\end{itemize}
\end{definition}

The main result of this paper is the following theorem, the proof is postponed to the last section.

\begin{theorem}\label{reducibilnost}
Let $G$ be a weak near-triangulation. If $G$ has  a block  $B$, such that $B$ is not outerplanar and the order of $B$ is different from 6, then $G$ is reducible. 
\end{theorem}

Note that it is not possible to extend the above theorem to weak near-triangulations that have only outerplaner blocks and  blocks on 6 vertices. We can see this by observing that 
 an octahedron embedds in the plane as a triangulation, and it is not dominated by one vertex, and therefore also not reducible. Another example is given in Fig.~\ref{slika} where the graph in case $(a)$ is not reducible. 
Observe also that every near-triangulation on $n$ vertices with domaination number $n/3$ is not reducible, 
examples of such outerplaner near-triangulations are exhibited in \cite{mathesontarjan}. In fact, Theorem \ref{reducibilnost} and Corollary~\ref{n3weak} imply the following. 

\begin{corollary} Every near-triangulation that attains the bound $n/3$ from Theorem \ref{mathesontarjan} 
 is  outerplanar or  has exactly 6 vertices.
\end{corollary}

\proof
 By Theorem \ref{reducibilnost}, every near-triangulation $G$ which is not outerplaner and has $n\neq 6$ vertices is reducible. 
So there exists a set $D$, such that $G-D$ is a weak near-triangulation, and $D$ is dominated by one vertex.  
The reduced graph $G-D$ has, according to Corollary \ref{n3weak}, a dominating set containing at most $(n-|D|)/3$ vertices. 
The result follows from $|D|\geq 4$. \qed

It is easy to prove that every near-triangulation $G$ on 6 vertices has a vertex $x$, such that the closed neighborhood of $x$ 
contains all internal vertices of $G$. The reader may prove this by a case analysis: the outer cylce of $G$ has 5, 4 or 3 vertices.

\begin{corollary}
Every plane triangulation on $n>6$ vertices has  a dominating set of size at most 17n/53.
\end{corollary}

\proof
Let $G$ be a plane triangulation on $n$ vertices. We apply Theorem \ref{reducibilnost} until the obtained graph $G'$ is irreducible. More precisely, $G'=G-(D_1\cup\ldots\cup D_k)$, 
where $D_i$ is a set of vertices in  $G-(D_1\cup\ldots\cup D_{i-1})$, such that there exists a vertex $x_i\in D_i$ with the property $D_i\subseteq N[x_i]$. Moreover, for all $i\leq k$, $|D_i|\geq 4$. 
By Theorem \ref{reducibilnost}, $G'$ is a weak near-triangulation, and 
every block of $G'$ is outerplanar or has exactly 6  vertices. If $G'$ contains all three external vertices of $G$, then $G'$ is a triangulation, and thus has exactly one block. 
If $G'$ is a block of order 6, then at least one reduction was done to get from $G$ to $G'$. The vertices of $V(G)\setminus V(G')$ are dominated by a dominating set of size at most 
$\frac 14 |V(G)\setminus V(G')|$. The block $G'$ of 
order 6 is dominated by 2 vertices. It follows that $G$ has a dominating set of size at most $3n/10$. If $G'$ is outerplaner, then it is a triangle. 
In this case $G$ has a dominating set 
of size at most $2n/7$. 

Assume therefore that at least one external vertex of $G$ is not in $G'$. We claim that every external vertex of $G'$ is adjacent to a vertex in $V(G)\setminus V(G')$. 
Let $x$ be an external vertex of $G'$. If $x$ is also an external vertex of $G$, the claim follows from the assumption. If $x$ is an internal vertex of $G$, then let $xy$ be an external edge of $G'$ incident to $x$. Since $x$ is an internal vertex of $G$, the edge $xy$ is an internal edge of $G$, and so there are two facial triangles in $G$ incident to $xy$, and since $xy$ is an external edge of $G'$, one of them is not a triangle of $G'$. So both vertices  $x$ and $y$ are adjacent to a vertex in $V(G)\setminus V(G')$. This proves the claim. 
Let $q=|V(G')|$.  By Corollary \ref{n3weak} $G$ has a dominating set of size at most $(n-q)/4+q/3$, this is one way to dominate $G$. Since every non-outerplanar block of $G'$ 
of order 6 has a vertex, 
that dominates all internal vertices of this block, we can choose a vertex in each block of order 6 to dominate internal vertices of $G'$. External vertices of $G'$ are, according to the claim above, dominated by vertices in $V(G)\setminus V(G')$. So we construct a dominating set by choosing all vertices in  $V(G)\setminus V(G')$ and a vertex in each block of order  $6$ to dominate internal vertices of $G'$. 
Since any two distinct blocks share at most one vertex, it follows that at most $|V(G')|/5$ vertices have been chosen in blocks of order 6. This construction gives a dominating set of size at most $n-q+ q/5$. For all $q$ we have $\min\{(n-q)/4+q/3,n-q+ q/5\}\leq  17n/53$.
\qed

\section{Operations on weak near-triangulations}

We start with definitions of a problematic and a bad vertex of a plane graph. 

\begin{definition}
Let $G$ be a plane graph and $u\in V(G)$. 
We say that $u$ is a {\em problematic vertex} in $G$ if $G-u$ is not a weak near-triangulation.  
\end{definition}

\begin{definition}\label{bad}
Let $G$ be a plane graph and $u\in V(G)$ a problematic vertex in $G$. 
We say that $x$ is a {\em $u$-bad vertex} in $G$ if $x$ is not contained in a triangle of $G-u$.  
\end{definition}
If $G$ is a weak near-triangulation, then every $u$-bad vertex in  $G$ is an external vertex of $G$, and it is adjacent to $u$.  Moreover, if $x$ is a $u$-bad vertex, then there exists a facial triangle $uxy$ in $G$. See Fig.\ref{sibkaskorajtriangulacija}, where $u$ is a problematic vertex in $G$, and there are two $u$-bad vertices in $G$ (the white vertices).
We skip the proof of the following lemma, it is similar to the proof of Lemma \ref{blok}.



\begin{lemma}\label{dajsklic}
Let $G$ be a weak near-triangulation and $x$ an external vertex of $G$. Every bounded face of $G-x$ is triangular. 
Moreover, if $y$ is adjacent to $x$ in $G$, then $y$ is an external vertex of  $G-x$. 
\end{lemma}
Lemma \ref{dajsklic} has an immediate corollary. 

\begin{corollary}\label{trikotnalica}
Let $G$ be weak near-triangulation, $D\subseteq V(G)$ a set of vertices such that $D\subseteq N[y]$  for  some $y\in D$. 
If $D$ contains at least one external vertex of $G$, then every bounded face of $G-D$ is triangular.  
\end{corollary}

\begin{lemma}\label{osnovna}
Let $G$ be a weak near-triangulation and $u$ a problematic external vertex in $G$. Let 
$T$ be the set of all $u$-bad verteces in $G$. Then $G-(\{u\}\cup T)$ is a weak near-triangulation.
\end{lemma}

\proof
Since every $u$-bad vertex in $G$ is adjacent to $u$, we find (by applying Corollary \ref{trikotnalica}) that every face of $G-(\{u\}\cup T)$ is triangular, except possibly the unbounded face. 
When removing vertices that are not contained in a triangle  of $G-u$ we obtain a graph in which every vertex is  contained in a triangle (note that this is possibly an empty graph).\qed

\noindent If $u$ is an internal vertex, then the following lemma applies. 

\begin{lemma}\label{osnovna1}
Let $G$ be a weak near-triangulation and $u$ a problematic internal vertex in $G$. Let 
$T$ be the set of all $u$-bad verteces in $G$, and suppose that $T\neq \emptyset$. Then $G-(\{u\}\cup T)$ is a weak near-triangulation.
\end{lemma}

\proof We note that every $u$-bad vertex in $G$ is an external vertex of $G$. The rest of the proof is the same as the proof of Lemma \ref{osnovna}. \qed

\begin{lemma}\label{lepljenje}
Let $G$ be a weak near-triangulation and $D\subseteq X\subseteq V(G)$ such that  
\begin{itemize}
\item[(i)] $G-X$ is a weak near-triangulation
\item[(ii)]  $ D\subseteq N[y]$ for some $y\in D$ 
\item[(iii)]  $D$ contains  at least one external vertex of $G$
\end{itemize}
then $G-D$ is a weak near-tirangulation if and only if every vertex of $X\setminus D$ is contained in a triangle of $G-D$. 

\end{lemma}

\proof
If a vertex of $X\setminus D$ is not contained in a triangle of $G-D$, then $G-D$ is not a WNT.

Suppose that every vertex in $X\setminus D$ is contained in a triangle of $G-D$. By (i) $G-X$ is a WNT, so every vertex of $G-X$ is contained in a triangle of $G-X$ and hence also in a triangle of $G-D$. It follows that every vertex of $G-D$ is contained in a triangle of $G-D$. We also see that from  (ii) and  (iii) together with Corollary \ref{trikotnalica} follows that every face of $G-D$ is triangular, except possibly the unbounded face. 
\qed

\begin{lemma} \label{trikotnik}
Let $G$ be a weak near-triangulation, and $X\subseteq V(G)$ a set  such that $G-X$ is a weak near-triangulation. Suppose that $uvw$ is a triangle of $G$, and let 
$Y$ be the set of vertices of $G$ contained in the interior of  $uvw$. Suppose that $X\cap Y=\emptyset$,  $u\in X$ and $v,w\notin X$.  
If $y\in Y$ is a problematic vertex in $G-X$, then every $y$-bad vertex in $G-X$ is adjacent to $u$ in $G$. 
\end{lemma}

\proof
Let $G,X,Y,y$ and $u,v,w$ be as declared in the lemma. Let $x$ be a $y$-bad vertex in $G-X$. 
If $x=v$ or $x=w$ there is nothing to prove, because $uvw$ is a triangle in $G$. Assume that $x\notin\{v,w\}$, and note that $x$ is adjacent to $y$, 
 and therefore lies  in the interior of $uvw$ in $G$.  
Since $x$ is a $y$-bad vertex in $G-X$, we find that there is a facial triangle $xyt$ in $G-X$.   
The edge $xt$ is an internal edge of $G$ and so there is a  facial triangle $xty'$ in $G$, such that $y'\neq y$. 
However, the triangle $xty'$ is not a triangle in $G-X$, for otherwise $x$ is not a $y$-bad vertex in $G-X$. 
It follows  that $y'$ is not in  $G-X$, and therefore $y'=u$ (because $y'$ is adjacent to $x$ which is an internal vertex of $uvw$).    
\qed

\section{Reducibility of weak near-triangulations} \label{miki}

In this section we prove Theorem \ref{reducibilnost}.
Let $G$ be a  weak near-triangulation and $B$  a block of $G$ that is not outerplanar. Suppose also that $B$ is not a block on 6 vertices.
Let $u$ be an internal vertex of $B$ so that $u$ has  at least two external neighbors (neighbors incident to the unbounded face of $G$). Let 
$u_1,\ldots,u_n$ be the external neighbors of $u$, and let $R_1,\ldots,R_n$ be  regions in the plane 
where the region $R_k$ is bounded by the outer cycle of $B$ and the edges $uu_k$ and $uu_{k+1}$ (indices are calculated  modulo $n$), if $n>2$ two consecutive regions share an edge and if $n=2$ they share two edges. We may choose $u$ so that at most one of the regions $R_k,k\leq n$ contains an internal vertex of $B$ different from $u$.  We note that this choice of 
$u$ only becomes relevant in subsection \ref{uniproblematicen}.

The proof is divided into four main cases (depending on the number of $u$-bad vertices in $G$ ) and several subcases. 
For each main case there is a subsection. In the proof we work with weak near-triangulation $G$, and with subgraphs of $G$. 
When we say "adjacent" we mean "adjacent in $G$" (as opposed to "adjacent in a subgraph of $G$"), and 
unless otherwise stated, 
a facial triangle means a facial triangle in $G$.

\subsection{There are at least three $u$-bad vertices in $G$ }

Suppose that $u$ has at least three $u$-bad vertices in $G$. Let $T$ be the set of all $u$-bad vertices in $G$, and define $D=\{u\}\cup T$. 
 By Lemma \ref{osnovna1},  $G-D$ is a WNT.  Since  $|D|\geq 4$ and $D\subseteq N[u]$,  $G$ is reducible.

\subsection{There are exactly two $u$-bad vertices in $G$ }

Suppose that $u$ has exactly two $u$-bad vertices in $G$, and let $x_1,x_2$ be $u$-bad vertices in $G$. By Lemma \ref{osnovna1}, $G-\{u,x_1,x_2\}$ is a WNT. 
Observe that $u$ is an internal vertex of $G$, so $ux_1$ and $ux_2$ are internal edges of $G$, and therefore there are two facial triangles $ux_1w_1$ and $ux_1z_1$ containing edge $ux_1$,  and two facial triangles $ux_2w_2$ and $ux_2z_2$ containing edge $ux_2$ (note that vertices $w_1,x_1,z_1$ are not necessarily distinct from $w_2,x_2,z_2$). Since $x_1$ and $x_2$ are $u$-bad vertices in $G$ neither $x_1$ nor $x_2$ is contained in a  triangle of $G-u$. It follows that $x_1w_1, x_1z_1$ and $x_2w_2,x_2z_2$ are external edges of $G$, and therefore  $w_1,w_2,x_1,x_2$ are external vertices of $G$. 
Note also that if $a\neq u$ is a common neighbor of $x_1$ and $x_2$ in $G$, then 
$a\in \{w_1,w_2,z_1,z_2\}$, because the facial triangles containing edges $ax_2$ and $ax_1$  also contain $u$ (because $x_1$ and $x_2$ are $u$-bad vertices in $G$). 
 We distinguish three possibilities. 1. $x_1$ and $x_2$ have a common neighbor different from $u$. 
2. $x_1$ and $x_2$ are  adjacent, and $u$ is their only common neighbor.
3.  $x_1$ and $x_2$ are not adjacent, and $u$ is their only common neighbor.

\subsubsection{$x_1$ and $x_2$ have a common neighbor different from $u$} Suppose that $x_1$ and $x_2$ have a common neighbor $t$ different from $u$. 
Observe that $t$ is adjacent to $u$, because $t\in \{w_1,w_2,z_1,z_2\}$. 
If $t$ is not problematic in  $G-\{u,x_1,x_2\}$, then let $D=\{u,x_1,x_2,t\}$. Since $G-D$ 
is a WNT and $D\subseteq N[t]$, $G$ is reducible. If $t$ is problematic in  $G-\{u,x_1,x_2\}$, then let $T$ be the set of all $t$-bad vertices in  $G-\{u,x_1,x_2\}$, and define 
$D=\{u,x_1,x_2,t\}\cup T$. By Lemma \ref{osnovna}, $G-D$ is a WNT. Since $D\subseteq N[t]$, $G$ is reducible. \label{skupensosed}

\subsubsection{$x_1$ and $x_2$ are adjacent, and  $u$ is their only common neighbor}
Observe that $ux_1x_2$ is a facial triangle in this case. Let $ux_1w_1$ and $ux_2w_2$ be facial triangles such that $w_1\neq x_2$ and $w_2\neq x_1$. 
Since  $u$ is the only common neighbor of $x_1$ and $x_2$ we have  $w_1\neq w_2$. 
If $w_1$ is not problematic in $G-\{u,x_1,x_2\}$, then define $D=\{u,x_1,x_2,w_1\}$. 
By the definition (of a problematic vertex) $G-D$ is a WNT. Since $|D|\geq 4$ and $ D\subseteq N[u]$, $G$ is reducible. If $w_2$ is not problematic in  $G-\{u,x_1,x_2\}$, the reduction is 
analogous as above, so assume that both $w_1$ and $w_2$  are problematic in  $G-\{u,x_1,x_2\}$.

Suppose that $w_1$ is a $w_2$-bad vertex in $G-\{u,x_1,x_2\}$, and that $w_2$ is a $w_1$-bad vertex in $G-\{u,x_1,x_2\}$.  In this case
$w_1$ and $w_2$ are adjacent, moreover $x_1$ and $x_2$ are in the exterior of triangle $w_1w_2u$, because $x_1$ and $x_2$ are external vertices of $G$. 
If  $w_1w_2x$ is in a facial triangle and $x\neq u$, 
then $x$  lies either in the interior of triangle $w_1w_2u$ or in the exterior. 
If $x$ is in the interior of $w_1w_2u$, then  the edge $xw_1$ is incident to a triangular face $xw_1t$. 
If $t\notin \{u,w_2\}$, then $xw_1t$ is a triangle in $G-\{u,x_1,x_2,w_2\}$, and so $w_1$ is not a $w_2$-bad vertex in $G-\{u,x_1,x_2\}$, contrary to the assumption.
It follows that $xw_1u$ and $xw_1w_2$ are facial triangles in $G$. Similarly, since $w_2$ is a $w_1$-bad vertex in $G-\{u,x_1,x_2\}$, we find that 
$xw_2u$ is a facial triangle in $G$.  
This in particular implies that $x$ is the only vertex of $G$  in the interior of triangle $w_1w_2u$.

\begin{figure}[htb] 
\begin{center}
\begin{tikzpicture}[scale=0.7]

\path node at (0,2) {$(a)$};
\filldraw (0,0) circle (2pt);
\filldraw (-1,1) circle (2pt);
\filldraw (1,1) circle (2pt);
\filldraw (-2,-2) circle (2pt);
\filldraw (2,-2) circle (2pt);


\filldraw (0,-3) circle (2pt);
\path node at (0,-2.5) {$y$};
\draw  (0,-3) to (-2,-2);
\draw  (0,-3) to (2,-2);

\draw  (0,0) to (-1,1);
\draw  (0,0) to (1,1);
\draw  (0,0) to (-2,-2);
\draw  (0,0) to (2,-2);
\draw  (-1,1) to (1,1);
\draw  (1,1) to (2,-2);
\draw  (-2,-2) to (-1,1);
\draw  (-2,-2) to (2,-2);

\path node at (0,0.5) {$u$};
\path node at (-1.5,1) {$x_1$};
\path node at (1.5,1) {$x_2$};
\path node at (-2.5,-2) {$w_1$};
\path node at (2.5,-2) {$w_2$};

\path node at (7,2) {$(b)$};
\filldraw (7,0) circle (2pt);
\filldraw (6,1) circle (2pt);
\filldraw (8,1) circle (2pt);
\filldraw (5,-2) circle (2pt);
\filldraw (9,-2) circle (2pt);

\filldraw (7,-1) circle (2pt);
\path node at (7,-1.5) {$x$};
\draw  (7,0) to (7,-1);
\draw  (7,-1) to (5,-2);
\draw  (7,-1) to (9,-2);


\draw  (7,0) to (6,1);
\draw  (7,0) to (8,1);
\draw  (7,0) to (5,-2);
\draw  (7,0) to (9,-2);
\draw  (6,1) to (8,1);
\draw  (8,1) to (9,-2);
\draw  (5,-2) to (6,1);
\draw  (5,-2) to (9,-2);

\path node at (7,0.5) {$u$};
\path node at (5.5,1) {$x_1$};
\path node at (8.5,1) {$x_2$};
\path node at (4.5,-2) {$w_1$};
\path node at (9.5,-2) {$w_2$};


\path node at (0,-4) {$(c)$};

\filldraw (0,-6) circle (2pt);
\filldraw (-1,-5) circle (2pt);
\filldraw (1,-5) circle (2pt);
\filldraw (-2,-8) circle (2pt);
\filldraw (2,-8) circle (2pt);



\draw  (0,-6) to (-1,-5);
\draw  (0,-6) to (1,-5);
\draw  (0,-6) to (-2,-8);
\draw  (0,-6) to (2,-8);
\draw  (-1,-5) to (1,-5);
\draw  (1,-5) to (2,-8);
\draw  (-2,-8) to (-1,-5);
\draw  (-2,-8) to (2,-8);

\path node at (0,-5.5) {$u$};
\path node at (-1.5,-5) {$x_1$};
\path node at (1.5,-5) {$x_2$};
\path node at (-2.5,-8) {$w_1$};
\path node at (2.5,-8) {$w_2$};


\path node at (7,-4) {$(d)$};

\filldraw (7,-6) circle (2pt);
\filldraw (6,-5) circle (2pt);
\filldraw (8,-5) circle (2pt);
\filldraw (5,-8) circle (2pt);
\filldraw (9,-8) circle (2pt);

\filldraw (7,-7) circle (2pt);
\path node at (7,-7.5) {$x$};
\draw  (7,-6) to (7,-7);
\draw  (7,-7) to (5,-8);
\draw  (7,-7) to (9,-8);

\filldraw (7,-9) circle (2pt);
\path node at (7,-8.5) {$y$};
\draw  (7,-9) to (5,-8);
\draw  (7,-9) to (9,-8);

\draw  (7,-6) to (6,-5);
\draw  (7,-6) to (8,-5);
\draw  (7,-6) to (5,-8);
\draw  (7,-6) to (9,-8);
\draw  (6,-5) to (8,-5);
\draw  (8,-5) to (9,-8);
\draw  (5,-8) to (6,-5);
\draw  (5,-8) to (9,-8);

\path node at (7,-5.5) {$u$};
\path node at (5.5,-5) {$x_1$};
\path node at (8.5,-5) {$x_2$};
\path node at (4.5,-8) {$w_1$};
\path node at (9.5,-8) {$w_2$};

\end{tikzpicture}
\caption{Cases $(a)$ to $(d)$.}\label{slika}
\end{center}
\end{figure}
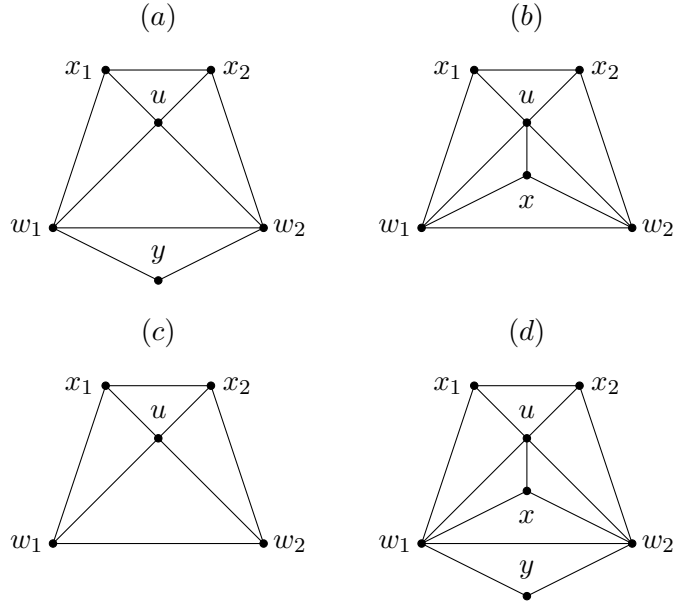

Suppose now that $w_1w_2y$ is a facial triangle such that $y$  is in the exterior of $w_1w_2u$. 
Since $w_1,w_2,x_1,x_2$ are external vertices of $G$, $y\neq x_1$ and $y\neq x_2$. 
If $yw_1$ is not an external edge of $G$, then there is a facial triangle $yw_1t$, 
with $t\neq x_1,x_2$ (note again that  $w_1,w_2,x_1,x_2$ are external vertices of $G$). This contradicts the assumption that $w_1$ is a $w_2$-bad vertex in $G-\{u,x_1,x_2\}$. 
Therefore $w_1y$ is an external edge of $G$, and similarly $w_2y$ is  an external edge of $G$. 
 Since $yw_1, yw_2, w_1x_1,x_1x_2$ and $x_2w_2$ are external edges of $G$, and also 
$w_1w_2$ is an external edge  if there is no facial triangle $w_1w_2y$ (with $y$ in the exterior of $w_1w_2u$), we find that $B$ is isomorphic to one of the four graphs  
shown in Fig.~\ref{slika}. $B$ is a block on 6 vertices in cases (a) and (b), which contradicts our assumptions.  $B$ is a block on 5  vertices in case (c), and since  $G-\{u,x_1,x_2\}$ is a WNT, $w_1$ 
is contained in a facial triangle $w_1ab$ of  $G-\{u,x_1,x_2\}$. Clearly, $a,b\notin B$, and therefore $w_1$ is not a $w_2$-bad vertex in  $G-\{u,x_1,x_2\}$ (contrary to the assumption). 
In case (d), $B$ is a block of order 7. If there is a  $w_1$-bad vertex $z\notin\{ x,w_2\}$ in $G-\{u,x_1,x_2\}$, then let $T$ be the set of all $w_1$-bad vertices in $G-\{u,x_1,x_2\}$, and define 
$D=(\{w_1,x_1\}\cup T)\setminus \{w_2\}$. Note that $x$ is a $w_1$-bad vertex in $G-\{u,x_1,x_2\}$ and so $x,z,x_1,w_1\in D$. It follows that
 $|D|\geq 4$, moreover $D\subseteq N[w_1]$. Since $G-\{u,x_1,x_2\}$ is a WNT, Lemma \ref{osnovna} implies that $G-(\{u,x_1,x_2,w_1\}\cup T)$ is a WNT. 
Since $(\{u,x_1,x_2,w_1\}\cup T)\setminus D=\{u,x_2,w_2\}$ we find, by applying Lemma \ref{lepljenje}, that $G-D$ is a WNT, and therefore $G$ is reducible.  (When applying 
 Lemma \ref{lepljenje} we set $X=\{u,x_1,x_2,w_1\}\cup T$, and we note that $X-D$ induces a triangle.)
If $x$ and $w_2$ are the only $w_1$-bad vertices in $G-\{u,x_1,x_2\}$, then define  $D=\{u,x_1,x_2,w_1,w_2,x\}$. Apply Lemma \ref{osnovna} to  $G-\{u,x_1,x_2\}$ (recall that 
$G-\{u,x_1,x_2\}$ is a WNT), and its 
 problematic vertex $w_1$. It follows that $G-D$ is a WNT, because $x$ and $w_2$ are the only $w_1$-bad vertices in $G-\{u,x_1,x_2\}$. Reducibility of $G$ follows from $D\subseteq N[u]$. 
This proves that  $w_1$ is not a $w_2$-bad vertex, or $w_2$ is not a $w_1$-bad vertex in  $G-\{u,x_1,x_2\}$.

 Assume, without loss of generality, that $w_2$ is not a $w_1$-bad vertex in  $G-\{u,x_1,x_2\}$. 
Let $T$ be the set of all $w_1$-bad vertices in $G-\{u,x_1,x_2\}$,  
and suppose first that $|T|\geq 2$. Define $D=\{x_1,w_1\}\cup T$. By Lemma \ref{osnovna}, $G-(\{u,x_1,x_2,w_1\}\cup T)$ is a WNT, and therefore by Lemma \ref{lepljenje}, 
$G-D$ is a WNT if and only if $u$ and $x_2$ are contained in a triangle of $G-D$. Since $ux_2w_2$ is a triangle in $G-D$, we conclude that $G-D$ is a WNT. 
Moreover $D\subseteq N[w_1]$ and $|D|\geq 4$, therefore $G$ is reducible.
Suppose now that there is exactly one  $w_1$-bad vertex  in $G-\{u,x_1,x_2\}$. Call it $z$, and note that $z\neq w_2$, 
and that   $G-\{u,x_1,x_2,w_1,z\}$ is a WNT by Lemma \ref{osnovna}.

The edge $w_1z$ is incident to a facial triangle in $G$ (and in $G-\{u,x_1,x_2\}$). 
Suppose that there exists a vertex $t\neq w_2$ such that $w_1zt$ is a facial triangle.
If $t$ is not problematic in   $G-\{u,x_1,x_2,w_1,z\}$ then let 
$D=\{x_1,w_1,z,t\}$. Since  $G-\{u,x_1,x_2,w_1,z,t\}$ is a WNT, and vertices $u$ and $x_2$ are contained in the 
triangle $ux_2w_2$ of $G-D$, we find by Lemma \ref{lepljenje}, that $G-D$ is a WNT. So $G$ is reducible, because $D\subseteq N[w_1]$. If $t$ is 
problematic in   $G-\{u,x_1,x_2,w_1,z\}$, then let $T$ be the set of all $t$-bad vertices in   $G-\{u,x_1,x_2,w_1,z\}$. 
Define $D=\{w_1,z,t\}\cup T$. By Lemma \ref{osnovna}, $G-\{u,x_1,x_2,w_1,z,t\}\cup T$ is a WNT, and so $G-D$ is a WNT by Lemma \ref{lepljenje}, due to the 
fact that vertices $u,x_1,x_2$ are contained in a triangle of $G-D$. Since $D\subseteq N[t]$ and $|D|\geq 4$, we find that $G$ is reducible. 
In both cases we found that $G$  is reducible, so assume that the only facial triangle incident to  $w_1z$ is the triangle $w_1zw_2$. We discuss two possibilities. 

First is when $z$ is in the interior of the triangle  $w_1w_2u$. Then $w_1z$ is an internal edge of $G$, so it's contained in two facial triangles. If one of them is $w_1zu$, then $z$ is adjacent to $u$ in $G$. 
In this case define $D=\{u,x_1,x_2,w_1,z\}$ and observe that $G-D$ is a WNT, and that $D\subseteq N[u]$, hence $G$ is reducible. Otherwise, if $w_1zu$ is not a facial triangle, there 
exists a vertex $t\neq w_2$ such that $w_1zt$ is a facial triangle, contradicting the assumption that $w_1zw_2$ is the only facial triangle containing the edge $w_1z$.

The second possibility is that $z$ is in the exterior of triangle $w_1w_2u$. 
Since $w_1zw_2$ is the only facial triangle incident to $w_1z$ we find that $w_1z$ is an external edge of $G$. Similarly, since $z$ is a 
$w_1$-bad vertex in  $G-\{u,x_1,x_2\}$, it is not contained in a triangle of $G-\{u,x_1,x_2,w_1\}$, and so $zw_2$ is an external edge of $G$. It 
follows that edges $w_2x_2,x_2x_1,x_1w_1,w_1z,zw_2$ are external edges of $G$, and so the union of these edges is the  boundary of  $B$.   
If there is no vertex in the interior of the triangle $w_1w_2u$, then $B$ is a block of order 6. 
Otherwise $u$ is adjacent to a vertex $y$ in the interior of $w_1w_2u$. 
 If $y$ is not problematic in $G-\{u,x_1,x_2\}$, then let $D=\{u,x_1,x_2,y\}$. 
$G-D$ is a WNT, by the definition of a problematic vertex, and since $D\subseteq N[u]$, $G$ is reducible. Assume therefore that $y$ is problematic in 
 $G-\{u,x_1,x_2\}$.
We apply Lemma \ref{trikotnik}, where $X=\{u,x_1,x_2\}$ and where the triangle $uw_1w_2$ takes the role of $uvw$ in Lemma \ref{trikotnik}. 
The lemma implies that every $y$-bad vertex in  $G-\{u,x_1,x_2\}$ is adjacent to $u$ in $G$. Let 
$T$ be the set of all $y$-bad vertices in  $G-\{u,x_1,x_2\}$. 
Define  $D=\{u,x_1,x_2,y\}\cup T$, and note that $D\subseteq N[u]$.  
 By Lemma \ref{osnovna},  $G-D$ is a WNT, therefore $G$ is reducible. 

\subsubsection{$x_1$ and $x_2$ are not adjacent, and $u$ is their only common neighbor}
 Let $w_1,z_1,w_2,z_2$ be pairwise distinct vertices such that 
 $ux_1w_1$ and $ux_1z_1$ are facial triangles containing $ux_1$,  and  $ux_2w_2$ and $ux_2z_2$ are facial triangles containing $ux_2$. 
Note that $x_1$ and $x_2$ are external vertices of $G$, and therefore $w_1$ is not adjacent to $z_1$.  Moreover, there is no path in $G$ that avoids $x_1,x_2$ and $u$ between 
$w_1$ and $z_1$ (again, because  $x_1$ and $x_2$ are external vertices).  Similarly  there is no path in $G$ that avoids $x_1,x_2$ and $u$ between $w_1$ and $w_2$, or 
there is no path in $G$ that avoids $x_1,x_2$ and $u$ between $w_1$ and $z_2$ (if both paths exist, then $x_1$ or $x_2$ is not an external vertex). Assume, without loss of generality, the latter. This, in particular, implies that $w_1$ is not adjacent 
to $z_2$ in $G$.

If $w_1$ is not problematic in  $G-\{u,x_1,x_2\}$, then we define $D=\{u,x_1,x_2,w_1\}$. Since $D\subseteq N[u]$ and $G-D$ is a WNT, $G$ is reducible. Assume therefore that $w_1$ is problematic in  $G-\{u,x_1,x_2\}$. 
If $w_1$ has at least two $w_1$-bad vertices in  $G-\{u,x_1,x_2\}$ then let 
$T$ be the set of $w_1$-bad vertices in  $G-\{u,x_1,x_2\}$ and define $D=\{w_1,x_1\}\cup T$.
Since $w_1$ is not adjacent to $z_2$, we find that $z_2\notin T$. 
It follows that $ux_2z_2$  is a triangle of $G-D$. We apply Lemma \ref{lepljenje} to  sets $X=\{u,x_1,x_2,w_1\}\cup T$ and $D$. Since 
$G-X$ is a WNT, and every vertex of $X\setminus D$ is contained in a triangle of $G-D$, we find that $G-D$ is a WNT. Since $D\subseteq N[w_1]$ and $|D|\geq 4$,  
$G$ is reducible. 
Assume now that there is exactly one  $w_1$-bad vertex in $G-\{u,x_1,x_2\}$, call this vertex $z$. 
Let $w_1zt$ be a facial triangle. If $t$ is not problematic in   $G-\{u,x_1,x_2,w_1,z\}$, then let $D=\{x_1,w_1,z,t\}$. 
Since  $z_2$ is not adjacent to  $w_1$, we find that both  $z$ and $t$ are distinct from  $z_2$.  
It follows that $u$ and $x_2$ are contained in a triangle of $G-D$, the triangle $ux_2z_2$. By Lemma \ref{lepljenje}, $G-D$ is a WNT, moreover $D\subseteq N[w_1]$. Therefore $G$ is reducible. The last possibility is that $t$ is problematic in  $G-\{u,x_1,x_2,w_1,z\}$. In this case let $T$ be the set of all $t$-bad vertices in  $G-\{u,x_1,x_2,w_1,z\}$, and define 
$D=\{t,z,w_1\}\cup T$. 
Since there is no path between $w_1$ and $z_2$ in $G$ that avoids $x_1,x_2$ and $u$, we find that  $z_2\notin T\cup\{z,t\}$. Analogous arguments prove that 
 $z_1\notin T\cup\{z,t\}$.  It follows that vertices $u,x_1,$ and $x_2$ are contained in triangles of $G-D$, these are triangles $ux_1z_1$ and $ux_2z_2$. 
By Lemma \ref{osnovna}, $G-(\{u,x_1,x_2,w_1,z,t\}\cup T)$ is a WNT. We apply Lemma \ref{lepljenje} to find that $G-D$ is a WNT. Moreover $D\subseteq N[t]$, so $G$ is reducible.  

\subsection{There is exactly one $u$-bad vertex in $G$ }

Suppose that  there is exactly one $u$-bad vertex in $G$, call it $x_1$. 
By Lemma \ref{osnovna1}, $G-\{u,x_1\}$ is a WNT. Let  $ux_1w_1$ and $ux_1z_1$ be facial triangles containing the edge $ux_1$.  
Note that $x_1z_1$ and $x_1w_1$ are external edges of $G$, and therefore $w_1$ and $z_1$ are external vertices of $G$.
If $w_1$ is problematic in $G-\{u,x_1\}$, then let $T$ be the set of all $w_1$-bad vertices in $G-\{u,x_1\}$. Define 
$D=\{u,x_1,w_1\}\cup T$. By Lemma \ref{osnovna}, 
$G-D$ is a WNT. Since $D\subseteq N[w_1]$, $G$ is reducible. 

Suppose that $w_1$ is not problematic in $G-\{u,x_1\}$, and suppose additionaly that 
$z_1$ is not problematic in $G-\{u,x_1,w_1\}$. Then define $D=\{u,x_1,w_1,z_1\}$, and observe that $G-D$ is a WNT. 
Since $D\subseteq N[u]$, $G$ is a reducible. 

It remains to prove that $G$ is reducible if $z_1$ is problematic in $G-\{u,x_1,w_1\}$. If $z_1$ is adjacent to $w_1$ in $G$, then let 
$T$ be the set of all $z_1$-bad vertices in $G-\{u,x_1,w_1\}$. Define $D=\{u,w_1,x_1,z_1\}\cup T$ and observe that $D\subseteq N[z_1]$. By Lemma \ref{osnovna}, 
$G-D$ is a WNT, and so $G$ is reducible. From now on we  assume that $w_1$ and $z_1$ are not adjacent in $G$.


If $\deg_B(w_1)=3$, then  $w_1$ is contained in a triangle  of   $G-\{u,x_1\}$, say $w_1ab$, such that $a,b\notin B$. In this case let $T$ be the set of all $z_1$-bad vertices in   $G-\{u,x_1,w_1\}$, and define  $D=\{u,x_1,z_1\}\cup T$. Since   $G-(\{u,x_1,w_1,z_1\}\cup T)$ is a WNT, also   $G-D$ is a WNT by Lemma \ref{lepljenje} (note that $a,b\notin T$, and 
so $w_1$ is contained in a triangle of $G-D$). Since $D\subseteq N[z_1]$, we find that $G$ is reducible. 
Therefore $\deg_B(w_1)\geq 4$.  

Now we can prove the following claim: every $z_1$-bad vertex in $G-\{u,x_1,w_1\}$ is adjacent to $w_1$ in $G$. If not, then the set of  
$z_1$-bad vertices in  $G-\{u,x_1,w_1\}$ that are not adjacent to $w_1$ in $G$ is nonempty, call this set $T$. Also, let $X$ be the set that contains $u,x_1,w_1,z_1$ and all $z_1$-bad vertices in $G-\{u,x_1,w_1\}$. By Lemma \ref{osnovna}, $G-X$ is a WNT.  Define $D=\{u,x_1,z_1\}\cup T$, and observe  that $G-D$ contains all 
neighbors of $w_1$ in $G$, except $x_1$ and $u$. Moreover, every vertex (distinct from $w_1$) in $X\setminus D$ is adjacent to $w_1$ in $G-D$. 
Since $\deg_B(w_1)\geq 4$, and $w_1x_1$ is an external edge of $B$, and $ux_1w_1$ is a facial triangle, we find that every vertex in $X\setminus D$ is contained  in a triangle 
of $G-D$. Hence, by Lemma \ref{lepljenje}, $G-D$ is a WNT. Since $D\subseteq N[z_1]$, $G$ is reducible. This proves the claim.

If every $z_1$-bad vertex in  $G-\{u,x_1,w_1\}$ is adjacent to $u$, then let $T$ be the set of all $z_1$-bad vertices in $G-\{u,x_1,w_1\}$. 
Define  $D=\{u,x_1,w_1,z_1\}\cup T$. Since $D\subseteq N[u]$, and $G-D$ is a WNT, $G$ is reducible.  Assume therefore that a $z_1$-bad vertex is not adjacent to $u$.

Let $y$ be a $z_1$-bad vertex in $G-\{u,x_1,w_1\}$ not adjacent to $u$ in $G$. Since $w_1z_1\notin E(G)$ and $yu\notin E(G)$ there exist exactly  one vertex adjacent to $y$ in the interior of 4-cycle $yw_1uz_1$  (for otherwise 
$y$ is not a $z_1$-bad vertex in  $G-\{u,x_1,w_1\}$). Call this vertex $z$ and observe that $zyw_1$ and $zyz_1$ are facial triangles. If $z$ and $u$ are not adjacent then  let $w_1zy_1$ be the facial triangle such that $y_1\neq y$. Let $T$ be the set of all $z_1$-bad vertices in $G-\{u,w_1,x_1\}$, and define $D=(\{u,x_1,z_1\}\cup T)\setminus \{z,y_1\}$. Since $G-\{u,w_1,x_1,z_1\}\cup T$ is a WNT, also $G-D$ is a WNT (by Lemma \ref{lepljenje}). 
As $y\in D$, we have $|D|\geq 4$ and so $G$ is reducible. It follows that $z$ and $u$ are adjacent. Moreover the triangle $w_1uz$ is a facial triangle 
(we do the same  reduction as above if $w_1uz$ is not a facial triangle).

Suppose that there is a vertex $x$  in the interior of triangle $uzz_1$, such that $x$ 
is adjacent to $u$ in $G$.  If $x$ is not problematic in  $G-\{u,x_1,w_1\}$, then let $D=\{u,x_1,w_1,x\}$ and observe that 
$G-D$ is a WNT. Since $D\subseteq  N[u]$, $G$ is reducible. 
Suppose that $x$ is problematic in    $G-\{u,x_1,w_1\}$. By Lemma \ref{trikotnik}, every $x$-bad vertex  in  $G-\{u,x_1,w_1\}$ is  adjacent to $u$. 
 Let $T$ be the set of all $x$-bad vertices in  $G-\{u,x_1,w_1\}$. Define 
$D= \{u,x_1,w_1,x\}\cup T$ and observe that $D\subseteq N[u]$. By Lemma \ref{osnovna}, $G-D$ is a WNT, and therefore $G$ is reducible. 
It follows that $uz_1z$ is a facial triangle of $G$, and therefore   $z$ is also a $z_1$-bad vertex in  $G-\{u,x_1,w_1\}$. 

 If the set $T$ of all $z_1$-bad vertices in  $G-\{u,x_1,w_1\}$ contains more than two vertices, then 
let $D=(\{u,x_1,z_1\}\cup T)\setminus \{z,y\}$. Since $zyw_1$ is a triangle in $G-D$ we find,  by applying Lemma \ref{lepljenje}, that $G-D$ is a WNT. Since $D\subseteq N[z_1]$, $G$ is reducible. 
Assume therefore that $z$ and $y$ are the only $z_1$-bad vertices in  $G-\{u,x_1,w_1\}$, and so $G-\{u,x_1,w_1,z_1,z,y\}$ is a WNT. If $w_1y$ is not an external edge of $G$, then define $D=\{u,x_1,z_1,z\}$. 
Since $w_1y$ is not an externl edge of $G$, $y$ and $w_1$ are contained in a triangle of $G-D$. Lemma \ref{lepljenje} implies that $G-D$ is a WNT.  Reducibility of $G$ follows from $D\subseteq N[z_1]$. Assume therefore that $yw_1$ is an external edge of $G$. If $yz_1$ is not an external edge of $G$, then define $D=\{u,x_1,w_1,z\}$. Since $y$ and $z_1$ are contained in a triangle of $G-D$ (note that there is a facial triangle $yz_1t$, where $t\neq z$, and $t\neq w_1$ because $z_1$ is not adjacent to $w_1$, and $t\neq x_1$ because $z_1$ is an external vertex of $G$), we find that $G-D$ is a WNT. Reducibility follows from $D\subseteq N[u]$. Hence also $z_1y$ is an external edge of $G$, and therefore all four edges $w_1x_1,x_1z_1,z_1y$ and $yw_1$ are external edges of $G$. So $B$ is a block of order 6.

\subsection{There are no $u$-bad vertices in $G$ }
\label{uniproblematicen}

Suppose that $u$ has no $u$-bad vertices in $G$. Let $u_1$ be any external neighbor of $u$ in  $G$. Since 
$u$ is not problematic in $G$,  $u_1$ is contained in a triangle of $G-u$.  Note first that $\deg_B(u_1)>2$, because $u_1$ is adjacent to $u$ in $B$, and $u$ is an internal vertex of $B$. 

We claim the following: if  $\deg_B(u_1)=3$ then $G$ is reducible or $u_1$ is contained in a triangle  $u_1ab$, where $a,b\notin B$. 
Suppose that  $\deg_B(u_1)=3$.  Let $z_1, z_2$ and $u$ be neighbors of $u_1$ in $B$, where $z_1$ and $z_2$ are external vertices of $B$.
Since $u_1$ is contained in a triangle of $G-u$ we find that $u_1$ is contained in a triangle $u_1ab$, where $a,b\notin B$, or $z_1$ is adjacent to $z_2$. 
Suppose that $z_1$ is adjacent to $z_2$ and that $u_1$ is not contained in a triangle $u_1ab$, where $a,b\notin B$. Then every block $B'\neq B$ of $G$ containing $u_1$ is a $K_2$,  and therefore  $u_1z_1z_2$ is the only triangle containing $u_1$ in $G-u$. It follows that $z_1$ and $z_2$ are the only vertices that are potentially not in a triangle of $G-\{u,u_1\}$.

Suppose that $\deg(u)\geq 4$ or $z_1z_2$ is an internal edge of $G$.  Then all vertices 
of $G-\{u,u_1\}$ are contained in a triangle of   $G-\{u,u_1\}$, and so (note that by Corollary \ref{trikotnalica} all bounded faces of  $G-\{u,u_1\}$ are triangular), 
$G-\{u,u_1\}$ is a WNT. In this case the reductions are defined as follows. 
If $z_1$ is problematic in $G-\{u,u_1\}$, then let $T$ be the set of all $z_1$-bad vertices in $G-\{u,u_1\}$, and define $D=\{u,u_1,z_1\}\cup T$. 
Since $G-D$ is a WNT and $D\subseteq N [z_1]$, $G$ is reducible. 
If $z_1$ is not problematic in 
$G-\{u,u_1\}$ and $z_2$ is not problematic in $G-\{u,u_1,z_1\}$, then let $D=\{u,u_1,z_1,z_2\}$. Since $D\subseteq N[u]$, $G$ is reducible. 
If $z_2$ is problematic in $G-\{u,u_1,z_1\}$, then let $T$ be the set of all $z_2$-bad vertices in $G-\{u,u_1,z_1\}$. Define $D=\{u,u_1,z_1,z_2\}\cup T$. Since $z_1$ is adjacent to $z_2$,  $D\subseteq N[z_2]$. Since  $G-D$ is a WNT, $G$ is reducible. 

It remains to define reductions when $z_1z_2$ is an external edge of $G$ and $\deg(u)=3$. In  this case $B$ is a $K_4$ (recall that $\deg_B(u_1)=3$, therefore $z_1u_1$ and $z_2u_1$ are external edges of $B$), and $G-u$ is a WNT. 
If there is a $z_1$-bad vertex $t\neq u_1$ in $G-u$, then   let $T$ be the set of all $z_1$-bad vertices in $G-u$, and define $D=\{u,z_1\}\cup T$. Since $u_1\in T$ we have  $|T|\geq 2$, and 
therefore $|D|\geq 4$. 
Since  $G-u$ is a WNT, we find (by applying Lemma \ref{osnovna}) that $G-D$ is a WNT. 
Reducibility follows from $D\subseteq N[z_1]$. Therefore $u_1$ is the only $z_1$-bad vertex in $G-u$, and similarly $u_1$ is also the only $z_2$-bad vertex in $G-u$. 
It follows that there exist blocks $B_1$ and $B_2$ of $G$, both different from $B$, and both near-triangulations, such that $z_1\in B_1$ and $z_2\in B_2$. 
So we may define $D=\{u,u_1,z_1,z_2\}$. Observe that $G-D$ is a WNT because none of the vertices $u_1,z_1$ and $z_2$ have a bad vertex in $G-u$ which is not contained in $B$. 
Since $D\subseteq N[u]$, $G$ is reducible.  
This proves the claim, and therefore it remains to define reductions in case if $\deg_B(u_1)\geq 4$, or  $\deg_B(u_1)=3$ and  $u_1$ is contained in a triangle  $u_1ab$, where $a,b\notin B$. 
Assume in the sequal that   $\deg_B(u_1)\geq 4$, or  $\deg_B(u_1)=3$ and  $u_1$ is contained in a triangle  $u_1ab$, where $a,b\notin B$, moreover assume that this is true for any external neighbor $u_1$ of $u$. This in particular implies that every vertex of $G-u$ is contained in a facial triangle of $G-u$ which is also a facial triangle of $G$ 
(or equvalently:  every vertex of $G-u$ is contained in a facial triangle of $G$ not containing $u$). 

\subsubsection{$u_1$ is not problematic in $G-u$}

Suppose that $u_1$ is not problematic in $G-u$, and therefore $G-\{u,u_1\}$ is a WNT. The edge $uu_1$ is incident to two triangular faces of $G$. So there are vertices $x$ and $y$ such that $uu_1x$ and $uu_1y$ are facial triangles in $G$. By the choice of $u$ (see Section \ref{miki}), at least one of $x$ and $y$ is an external vertex of $B$ (and $G$), for otherwise more than one region $R_k, k\leq n$ 
contains an internal vertex of $B$. Assume that $x$ is an external vertex.  If $u_1y$ is an external edge of $G$, then let $u_2=y$, otherwise let $u_2=x$. 

We claim that $u_1$ is contained in a triangle of $G-\{u,u_2\}$.
If $\deg_B(u_1)=3$ then $u_1$ is contained in a triangle  $u_1ab$, where $a,b\notin B$ (in which case the claim is true). Assume now that $\deg_B(u_1)\geq 4$. Now if $u_2=y$ then $u_1u_2$ is an external edge of $B$, 
and $uu_1u_2$ is a facial triangle. Since $\deg_B(u_1)\geq 4$ and $B$ is a near-triangulation we find that $u_1$ is incident to at least three facial triangles of $B$. 
When we remove $u$ and $u_1$ from $G$ at least one facial triangle containing $u_1$ remains. If $u_2=x$ and  $u_1x$ is an external edge, then the argument is the same as above. So assume that $u_1x$ is also not an external edge (we already know that $u_1y$ is not an external edge). Now in this case $\deg_B(u_1)\geq 5$. Removing vertices $u$ and $u_2$ (which form a facial triangle with $u_1$) keeps at least one facial triangle containing $u_1$.  This proves the claim. 

If $u_2$ is a problamatic vertex in $G-\{u,u_1\}$, then let $T$ be the set of all $u_2$-bad vertices in  $G-\{u,u_1\}$. Since  $G-\{u,u_1\}$ is a WNT, it follows from 
Lemma \ref{osnovna} that  $G-(\{u,u_1,u_2\}\cup T)$ is a WNT. We define $D=\{u,u_1,u_2\}\cup T$.  Since $D\subseteq N[u_2]$, $|D|\geq 4$ and $G-D$ is a WNT, we find that 
$G$ is reducible. 

Assume that $u_2$ is not problematic in  $G-\{u,u_1\}$ and so $G-\{u,u_1,u_2\}$ is a WNT. Let $t\neq u_1$ be such that $uu_2t$ is a facial triangle. 
If $t$ is not problematic in $G-\{u,u_1,u_2\}$ then let $D=\{u,u_1,u_2,t\}$. We have $D\subseteq N[u]$ and $G-D$ is a WNT, so $G$ is reducible. So assume that $t$ is problematic 
in  $G-\{u,u_1,u_2\}$.

Let $T$ be the set of all $t$-bad vertices in  $G-\{u,u_1,u_2\}$. If $t$ and $u_1$ are adjacent in $G$ then let $D=\{u,u_1,u_2,t\}\cup T$. We find that 
$G-D$ is a WNT according to Lemma \ref{osnovna}. Moreover $D\subseteq N[t]$, and so $G$ is reducible. Assume therefore that $t$ is not adjacent to $u_1$. 
If no vertex of $T$ is adjacent to $u_1$ then let  $D=\{u,u_2,t\}\cup T$. It follows from Lemma \ref{lepljenje} and the fact that $u_1$ is contained in a triangle of $G-\{u,u_2\}$ (and so also in a triangle of $G-D$) that $G-D$ is a WNT. Since $D\subseteq N[t]$, $G$ is reducible.  Assume from now on that a vertex of $T$ is adjacent to $u_1$, and call this vertex $t_0$.

{\bf Case $(a)$.} Suppose that $t_0$ is not adjacent to $u$ in $G$. Then the interior of 4-cycle $t_0u_1ut$ contains exactly one vertex  adjacent to $t_0$, call it $z$
(if the interior of $t_0u_1ut$ contains two vertices adjacent to $t_0$, then $t_0$ is contained in a triangle of $G-\{u,u_1,u_2,t\}$, and hence $t_0$ is 
 not a $t$-bad vertex in  $G-\{u,u_1,u_2\}$).  
If $z$ and $u$ are not adjacent in $G$, then $u_1$ is contained in a triangle  $u_1zz'$, where $z'$ lies in the interior of  $u_1utz$. 
Let $T$ be the set of all $t$-bad vertices in $G-\{u,u_1,u_2\}$, and define $D=\{u,u_2,t\}\cup T\setminus\{z,z'\}$. Since $G-(\{u,u_1,u_2,t\}\cup T)$ is a WNT, we find (by applying Lemma \ref{lepljenje}) that 
$G-D$ is a WNT. Since $t_0\in D$, we have $|D|\geq 4$. Moreover$D\subseteq N[t]$, and therefore $G$ is reducible.

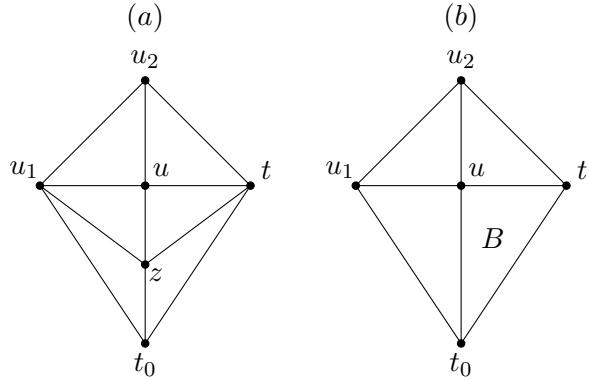
\begin{figure}[htb!] 
\begin{center}
\begin{tikzpicture}[scale=0.7]

\path node at (0,3.2) {$(a)$};
\filldraw (0,0) circle (2pt);
\filldraw (0,2) circle (2pt);
\filldraw (2,0) circle (2pt);
\filldraw (-2,0) circle (2pt);
\filldraw (0,-3) circle (2pt);
\filldraw (0,-1.5) circle (2pt);


\draw  (0,0) to (0,2);
\draw  (0,0) to (2,0);

\draw  (0,0) to (-2,0);
\draw  (-2,0) to (0,2);
\draw  (0,2) to (2,0);
\draw  (0,0) to (0,-3);
\draw  (-2,0) to (0,-3);
\draw  (2,0) to (0,-3);
\draw  (0,-1.5) to (2,0);
\draw  (0,-1.5) to (-2,0);

\path node at (0.3,0.3) {$u$};
\path node at (-2.3,0.3) {$u_1$};
\path node at (0,2.4) {$u_2$};
\path node at (2.3,0.3) {$t$};
\path node at (0,-3.4) {$t_0$};
\path node at (0.2,-1.7) {$z$};



\path node at (6,3.2) {$(b)$};
\filldraw (6,0) circle (2pt);
\filldraw (6,2) circle (2pt);
\filldraw (8,0) circle (2pt);
\filldraw (4,0) circle (2pt);
\filldraw (6,-3) circle (2pt);


\draw  (6,0) to (6,2);
\draw  (6,0) to (8,0);

\draw  (6,0) to (4,0);
\draw  (4,0) to (6,2);
\draw  (6,2) to (8,0);
\draw  (6,0) to (6,-3);
\draw  (4,0) to (6,-3);
\draw  (8,0) to (6,-3);

\path node at (6.3,0.3) {$u$};
\path node at (3.7,0.3) {$u_1$};
\path node at (6,2.4) {$u_2$};
\path node at (8.3,0.3) {$t$};
\path node at (6,-3.4) {$t_0$};
\path node at (6.6,-1) {$B$};

\end{tikzpicture}
\caption{Case $(a)$: $t_0$ is not adjacent to $u$,  Case $(b)$: $t_0$ is adjacent to $u$.}\label{primeri}
\end{center}
\end{figure}

Assume therefore that $z$ and $u$ are adjacent in $G$ (see  Fig.~\ref{primeri}, Case $(a)$). If the interior of the triangle $u_1uz$ contains a vertex, then  $u_1$ is contained in a triangle $u_1zz'$, where $z'$ lies in the interior of  $u_1uz$. In this case we can do the same reduction as in the previous case, by defining $D=\{u,u_2,t\}\cup T\setminus\{z,z'\}$. 
We may therefore assume that there is no vertex of $G$ in the interior of $uu_1z$.

 Suppose that $y$ is a vertex of $G$ in the interior of $utz$ adjacent to $u$. If
 $y$ is not problematic in $G-\{u,u_1,u_2\}$, we define $D=\{u,u_1,u_2,y\}$. Since $G-D$ is a WNT, and $D\subseteq N[u]$, $G$ is reducible. 
Otherwise, if $y$ is problematic in  $G-\{u,u_1,u_2\}$, then by Lemma  \ref{trikotnik}, every $y$-bad vertex in $G-\{u,u_1,u_2\}$ is adjacent to $u$ in $G$. 
In this case let $T$ be the set of all $y$-bad vertices in $G-\{u,u_1,u_2\}$ and define $D=\{u,u_1,u_2,y\}\cup T$. Since $G-D$ is a WNT, and $D\subseteq N[u]$, we find that $G$ is reducible. 
This proves that there is no vertex of $G$ in the interior of $utz$.  Observe that we proved that all triangles on  Fig.~\ref{primeri}, Case $(a)$, are facial triangles. 

 Let $T$ be the set of all $t$-bad vertices in $G-\{u,u_1,u_2\}$, and let  $X=\{u,u_1,u_2,t\}\cup T$.
If $|T|\geq 3$, then $D=X\setminus \{u_1,z,t_0\}$. Since $G-X$ is a WNT, also $G-D$ is a WNT, according to Lemma \ref{lepljenje}. Since $D\subseteq N[t]$, and $D\geq 4$, $G$ is reducible. 
Assume therefore that $z$ and $t_0$ are the only $t$-bad vertices in  $G-\{u,u_1,u_2\}$, and so $X=\{ u,u_1,u_2,t,z,t_0\}$. 

If $u_1t_0$ is not an external edge of $G$, there exists a facial triangle $u_1t_0x$, with $x\neq z$, $x\neq t$ (recall that $u_1$ is not adjacent to $t$) and $x\neq u_2$ (recall that $u_1$ is an external vertex). In this case we define $D=X\setminus \{u_1,t_0\}=\{u_2,u,t,z\}$. Since every vertex of $X-D$ is contained in a triangle of $G-D$ we find, by applying Lemma \ref{lepljenje}, that
 $G-D$ is a WNT. Since $D\subseteq N[u]$, $G$ is reducible. It follows that $u_1t_0$ is an external edge. 
With an analogous argument we prove that $u_1u_2$ is an external edge. If $u_2t$ is not an external edge of 
$G$, then there is a facial  triangle $u_2tx$, with $x\neq u$ and $x\neq u_1$ (recall that $u_1$ is not adjacent to $t$). If $x=t_0$, then $u_2$ and $t_0$ are adjacent in $G$.
If $u_2t_0$ is an external  edge of $G$, then all three edges $u_1u_2,u_2t_0$ and $t_0u_1$ are external edges of $G$, and therefore $B$ is a block of order 6 in $G$. 
Otherwise, if 
$u_2t_0$ is an internal edge, then there exists an $x\notin \{t, u_1\}$ such that $u_2t_0x$ is a facial triangle. In this case we define $D=X\setminus \{u_2,t_0\}$. Since $G-X$ is
 a WNT, also $G-D$ is a WNT. Moreover $D\subseteq N[u]$, so $G$ is reducible. If $x\neq t_0$, then define $D=\{u,u_1,z,t_0\}$. Since $u_2$ and $t$ are contained in a triangle of 
$G-D$, we find that $G-D$ is a WNT. Since $D\subseteq N[u_1]$, we find that $G$ is reducible. Therefore we may assume that $u_2t$ is also an external edge of $G$. 
Finally we argue that $t_0t$ is an external edge of $G$. If not, then there is a facial triangle $tt_0x$, with $x\neq z$. If $x=u_2$, then $u_2t$ is not an external edge,  contradicting the above assumption. Otherwise $x\neq u_2$ and $x\neq u_1$. 
So we define $D=X\setminus \{t,t_0\}$ and argue that $G-D$ is a WNT by refering to Lemma \ref{lepljenje}. Altogether, we proved that 
$u_1u_2,u_2t,tt_0,t_0u_1$ are external edges of $G$, therefore $B$ contains $u,u_1,u_2,t,z,t_0$ and no other vertices, hence $B$ is a block of order 6.

{\bf Case $(b)$.} Now we discuss the case when $t_0$ and $u$ are adjacent in $G$ (see Fig.~\ref{primeri}, Case $(b)$). If there is a vertex in the interior of $u_1ut_0$, then there is exactly one such vertex, for otherwise either 
$t_0$ is contained in a triangle of $G-\{u,u_1,u_2,t\}$ (and so $t_0$ is not a $t$-bad vertex) or there is a triangle $u_1zz'$, where $z$ and $z'$ are in the interior of $u_1ut_0$ 
(note that in this case $z$ and $z'$ are not adjacent to $t$).
If the latter happens, then let  $T$ be the set of all $t$-bad vertices in $G-\{u,u_1,u_2\}$ and define $D=(\{u,u_2,t\}\cup T)$. 
Since $u_1zz'$ is a triangle in $G-D$, we find that $G-D$ is a WNT (by Lemma \ref{lepljenje}). Since $D\subseteq N[t]$, we find that  $G$ is reducible. 
If there is exactly one vertex in the interior of $u_1ut_0$, then this vertex is a $u_1$-bad vertex in $G-u$, and so $u_1$ is problematic in $G-u$. This case is treated in Section \ref{u1jeproblematicen}. 
Assume therefore that there is no vertex  in the interior of $u_1ut_0$, so we have Case $(b)$ of Fig.~\ref{primeri}, where $B$ potentially contains some vertices of $G$.

Suppose that $t_0$ is the only $t$-bad vertex in $G-\{u,u_1,u_2\}$. In this case let $D=\{u,u_1,u_2,t,t_0\}$.
By Lemma \ref{osnovna}, $G-D$ is a WNT. Since $D\subseteq N[u]$, $G$ is reducible. Assume therefore that there are at least two  $t$-bad vertices in $G-\{u,u_1,u_2\}$. %

If $t_0u_1$ is not an external edge of $G$, then there is a facial triangle $t_0u_1x$, with $x\neq u$ and $x\neq t$ (recall that $u_1$ is not adjacent to $t$) and $x\neq u_2$ (recall that $u_1$ is an external vertex). If $x$ is a $t$-bad vertex in $G-\{u,u_1,u_2\}$  
 then we get case $(a)$ of Fig.~\ref{primeri}. Assume therefore that $x$ is not a $t$-bad vertex  in $G-\{u,u_1,u_2\}$. Let $T$ be the set of all 
 $t$-bad vertices in $G-\{u,u_1,u_2\}$, $X=\{u,u_1,u_2,t\}\cup T$ and $D=X\setminus \{u_1,t_0\}$.  By Lemma \ref{osnovna} we find that $G-X$ is a WNT, and therefore, by applying Lemma \ref{lepljenje} we see that $G-D$ is a WNT. Reducibility of $G$ follows from $D\subseteq N[t]$ and $|D|\geq 4$. We conclude that $u_1t_0$ is an external edge of $G$. 

If $u_1u_2$ is not an external edge of $G$, then there is a facial triangle $u_1u_2x$, with $x\neq u$, $x\neq t$ and $x\neq t_0$. Moreover $x$ is not adjacent to $t$ because $u_1$ and $u_2$ are external vertices of $G$. In this case we define the reduction 
by $D=X\setminus \{u_1,u_2\}$. This proves that both edges $t_0u_1$ and $u_1u_2$ are external edges of $G$. 
It follows that $\deg_B(u_1)=3$, and so $u_1$ is contained in a triangle $u_1ab$, with $a,b\notin B$ 
(see the assumption at the end of Section \ref{uniproblematicen}).  
Observe that $a$ and $b$ are not adjacent to $t$ (because $t\in B$). 
We define the reduction by  $D=X\setminus \{u_1\}$. Since $u_1$ is contained in a triangle of $G-D$ we see that $G-D$ is a WNT, and so the reduction is well defined. 

\subsubsection{$u_1$ is problematic in $G-u$.}\label{u1jeproblematicen}

 Note that there are no $u$-bad vertices in $G$, and so every 
vertex of $G-u$ is contained in a triangle of $G-u$.

Suppose that $u_1$ is problematic in $G-u$.
If there are at least two $u_1$-bad vertices in $G-u$, then let $T$ be the set of all $u_1$-bad vertices in $G-u$, 
and define   $D=\{u,u_1\}\cup T$.  
Since  every 
vertex of $G-u$ is contained in a triangle of $G-u$, we find that every $u_1$-bad vertex in $G-u$ is adjacent to $u_1$. 
By Corollary \ref{trikotnalica} every bounded face of $G-D$ is triangular. 
Moreover, removing $u_1$ and all $u_1$-bad vertices from $G-u$ produces a graph in which all vertices are contained in a triangle. 
Therefore $G-D$ is a WNT, and since  $D\subseteq N[u_1]$ and $|D|\geq 4$, $G$ is reducible. 
Assume therefore that there is exactly one $u_1$-bad vertex $x$ in $G-u$, and observe that the same arguments as above prove that $G-\{u,u_1,x\}$ is a WNT. Since every vertex of $G-u$ (and in particular $x$) is contained in a facial triangle of $G$ not containing $u$ 
(see the assumption at the end of Section \ref{uniproblematicen}), there exists a facial triangle  $u_1xt$ in $G$, where $t\neq u$.
(since $x$ is a $u_1$-bad vertex in $G-u$,  both vertices $u_1$ and $x$ are contained in the same facial triangle, and removing $u_1$ breaks this facial triangle). 

 First let us prove that 
$u$ and $t$ are not adjacent. If they are adjacent,  then let $T$ be the (possibly empty) set of all $t$-bad vertices in $G-\{u,u_1,x\}$ and define $D=\{u,u_1,x,t\}\cup T$. Since $D\subseteq N[t]$ and $G-D$ is a WNT we find that $G$ is reducible. This proves that $u$ is not adjacent to $t$.
 If $t$ is not problematic in 
$G-\{u,u_1,x\}$, then let $D=\{u,u_1,x,t\}$. By the defintion (of a problematic vertex) $G-D$ is a WNT, moreover $D\subseteq N[u_1]$, so $G$ is reducible.
Assume therefore that $t$ is problematic in $G-\{u,u_1,x\}$. First we will assume that there is a $t$-bad vertex $t_0$ in $G-\{u,u_1,x\}$ adjacent to $u$.

\bigskip

\noindent {\bf Case 1: A $t$-bad vertex is adjacent to $u$.}

\bigskip
\noindent There is at most one vertex in the interior of the 4-cycle $uu_1tt_0$ adjacent to $t_0$. If not, $t_0$ is contained in a triangle of $G-\{u,u_1,x,t\}$, and so $t_0$ is not a $t$-bad vertex in $G-\{u,u_1,x\}$. Similarly, 
there is at most one vertex in the interior of the 4-cycle $uu_1tt_0$ adjacent to $u$. If not, $u$ is contained in a triangle $uw_1w_2$, where both $w_1$ and $w_2$ are in the interior of  
 $uu_1tt_0$. In this case let $T$ be the set of all $t$-bad vertices in  $G-\{u,u_1,x\}$ and define $D=(\{u_1,x,t\}\cup T)\setminus \{w_1,w_2\}$.
Since $\{u_1,x,t,t_0\}\in D$ we have $|D|\geq 4$, and since $G-(\{u,u_1,x,t\}\cup T)$ is a WNT, we find (by applying Lemma \ref{lepljenje}) that $G-D$ is a WNT. Therefore assume that 
either there is no vertex of $G$ in the interior of  the 4-cycle  $uu_1xt$ adjacent to $t_0$ and $u$, or there is exactly one vertex in the interior of $uu_1xt$ adjacent to $t_0$ and $u$ (and note that this must be the same vertex). Call this vertex (if it exists) $z$,  and note that $u_1zu$, $uzt_0$ and $tzt_0$ are facial triangles in this case.

We claim that $t_0$ is the only $t$-bad vertex in $G-\{u,u_1,x\}$. Let us first exclude the possibility $|T|\geq 3$, where $T$ is the set of all $t$-bad vertices in $G-\{u,u_1,x\}$. 
The edge $uu_1$ is an internal edge of $G$, so there is a facial triangle $uu_1w$, with $w\neq t$. 
If $|T|\geq 3$, let $D=(\{x,t\}\cup T)\setminus\{w\}$. $|D|\geq 4$ and since $G-(\{u,u_1,x,t\}\cup T)$ is a WNT we find (by 
Lemma \ref{lepljenje}) that   $G-D$ is a WNT. Since $D\subseteq N[t]$, $G$ is reducible.   This leaves us with the possibility $|T|=2$, to be excluded next. 

Let $t_0$
 and $y$ be the only $t$-bad vertices in $G-\{u,u_1,x\}$. If there is no vertex in the interior of $uu_1tt_0$, then $u_1$ and $t_0$ are adjacent (recall that $u$ and $t$ are not adjacent). 
Let $uu_1w$ be the facial triangle, with $w\neq t_0$. If $w=y$, then let $D=\{u,u_1,x,t,t_0,y\}$. Since $D\subseteq N[u_1]$, $G$ is reducible. 
If $w\neq y$, the reduction is given by  $D=\{x,t,t_0,y\}$. Since $w\neq t$ (recal that $t$ and $u$ are not adjacent), and $w\neq x$ (recall that $u_1$ 
is an external vertex of $G$) we find that $uu_1w$ is a triangle in $G-D$, and so $G-D$ is a WNT (by Lemma \ref{lepljenje}). Since $D\subseteq N[t]$, $G$ is reducible. 
 Suppose now that there is a
vertex $z$ in the interior of $uu_1tt_0$ adjacent to $t_0, u,u_1$ and $t$. 
First assume that $z=y$. If $u_1$ and $t_0$ are adjacent, then define   $D=\{u,u_1,x,t,t_0,y\}$. We 
have $D\subseteq N[u_1]$ and $G-D$ is a WNT, so $G$ is reducible. Now assume that  $u_1$ and $t_0$ are not adjacent. In this case there is a facial triangle $uu_1w$, with $w\neq t_0$ and $w\neq y=z$. Also $w\neq x$ because $u_1$ is an external vertex of $G$, and $w\neq t$ since $u$ and $t$ are not adjacent.  In this case we define 
$D=\{t,x,y,t_0\}$. Observe that $D\subseteq N[t]$ and that $uu_1w$ is a triangle in $G-D$, so $G-D$ is a WNT (by Lemma  \ref{lepljenje}). This proves that 
 $G$ is reducible. 
It remains to check the case $z\neq y$. In this case we have 
 $D=\{x,t,t_0,y\}$, and so $uu_1z$ is a triangle in $G-D$. It follows that $G-D$ is a WNT (by Lemma \ref{lepljenje}). Reducibility of $G$ follows from $D\subseteq N[t]$. 
This proves that $t_0$ is the only $t$-bad vertex in $G-\{u,u_1,x\}$, as claimed.

 If $u_1$ and $t_0$ are adjacent then  $D=\{u,u_1,x,t,t_0\}$. Since $D\subseteq N[u_1]$, and $G-D$ is a WNT, $G$ is reducible. Suppose therefore that $u_1$ and $t_0$ are not adjacent. 
In this case there is a vertex $z$ is the interior of $uu_1tt_0$ adjacent to all four vertices of this 4-cycle.
Since $ut_0$ is an internal edge, there is a facial triangle $ut_0w$, with $w\neq z$. We have $w\neq t$ (because $u$ is not adjacent to $t$), and $w\neq u_1$  (because $u_1$ and $t_0$ are not adjacent by the assumption) and $w\neq x$ (for otherwise  $x$ is contained in the triangle $xtt_0$ in $G-\{u,u_1\}$, and so $x$ is not a $u_1$-bad vertex in $G-u$). If $z$ is not problematic in $G-\{u,u_1,x,t,t_0\}$ then $G-\{u,u_1,x,t,t_0,z\}$ is a WNT. Define $D=\{t,x,z,u_1\}$ and observe that $D\subseteq N[t]$. Since $ut_0w$ is a triangle in $G-D$ we find (by Lemma \ref{lepljenje}) that $G-D$ is a WNT. 
It remains to check the case when $z$ is   problematic in $G-\{u,u_1,x,t,t_0\}$. Let $T$ be the set of all $z$-bad vertices in  $G-\{u,u_1,x,t,t_0\}$. Define $D=\{z,u,t_0\}\cup T$. 
Since $G-\{u,u_1,x,t,t_0,z\}\cup T$ is a WNT, we find that $G-D$ is a WNT, because $u_1xt$ is a triangle in $G-D$. $G$ is reducible because $D\subseteq N[z]$.

\bigskip

\noindent {\bf Case 2: No $t$-bad vertex is adjacent to $u$.}

\bigskip
\noindent
Assume now that no $t$-bad vertex is adjacent to $u$. 
We will prove that $x$ and $u$ are not adjacent in $G$. For the purpose of a contradiction assume that $x$ and $u$ are adjacent. Then there are two possible drawings of the graph induced by $u,u_1,x,t$: either the triangle $u_1xt$ lies in the interior of the triangle $u_1xu$, or it lies in the exterior. 
If the triangle $u_1xt$ lies in the interior of the triangle $u_1xu$ then the edge $xt$
 is incident to a triangular face contained in the interior of the 4-cycle $u_1txu$. Let $xtw$ be the triangle bounding this triangular face. Since there are no multiple edges in $G$ and $u$ is not adjacent to $t$, we find that $w\notin \{u,u_1\}$. It follows that $x$ is contained in a  traingle of $G-\{u,u_1\}$, a contradiction (recall that $x$ is a $u_1$-bad vertex in $G-u$). 

Let's consider the case when  triangle $u_1xt$ lies in the exterior of the triangle $u_1xu$. We claim that,  the interior of $uu_1x$ contains at most one vertex of $G$. To see this let $w$ be a vertex in the interior of $uu_1x$ adjacent to $x$. The edge $xw$ is an internal edge, so it is incident to two facial triangles $xwy$ and $xwy'$. If $y$ or $y'$ is not $u$ or $u_1$, then $x$ 
is contained in a triangle of $G-\{u,u_1\}$, and so $x$ is not a $u_1$-bad vertex in $G-u$, a contradiction. Therefore $y=u_1$ and $y'=u$. If there is a vertex of $G$ in the interior of the tirangle $uu_1w$, then $u$ is contained in a triangle of $G-(\{u_1,x,t\}\cup T)$, where  $T$ is  the set of all $t$-bad vertices in   $G-\{u,u_1,x\}$. Define $D=\{u_1,x,t\}\cup T$. 
By Lemma \ref{osnovna},  $G-(\{u,u_1,x,t\}\cup T)$ is a WNT, and therefore $G-D$ is a WNT by Lemma \ref{lepljenje}. This proves that there is no vertex in the interior of $uu_1w$, and so $w$ is the only vertex in the interior of triangle $uu_1x$, as claimed.

We now discuss both cases: there is a vertex in the interior of $uu_1x$, or there is no such vertex. 
Suppose that $w$ is (the only) vertex in the interior of triangle $uu_1x$. Then $w$ is a $u_1$-bad vertex in $G-u$, this contradicts that $x$ is the only $u_1$-bad vertex in $G-u$ 
(see assumptions in the second paragraph  of  Section \ref{u1jeproblematicen}). It 
remains to check the case when there is no vertex in the interior of $uu_1x$. In this case $uu_1x$ is a facial triangle, and if $\deg_G(u)\geq 4$, then $u$ is contained in a triangle of 
 $G-(\{u_1,x,t\}\cup T)$, because $t$ is not adjacent to $u$, and no vertex in $T$ is adjacent to $u$ (and because $u$ is an internal vertex of $G$). In this case we define $D=\{u_1,x,t\}\cup T$. Since  $G-(\{u,u_1,x,t\}\cup T)$ is a WNT, we find  by an application of  Lemma \ref{lepljenje} 
that $G-D$ is a WNT. Since $D\subseteq N[t]$, $G$ is reducible. It follows that $\deg_G(u)=3$, and let $w\neq u_1,x$ be the third neighbor of $u$ in $G$. Clearly, $w$ is adjacent to $x$ 
 and $u_1$, because $u$ is an internal vertex and $G$ is a WNT. Let $T$ be the set of all $w$-bad vertices in $G-\{u,u_1,x\}$ (if any), and define $D=\{u,u_1,x,w\}\cup T$. By Lemma \ref{osnovna}, $G-D$ is a WNT, and since $D\subseteq N[w]$, we conclude that $G$ is reducible. This proves that $u$ and $x$ are not adjacent in $G$. 

Since no $t$-bad vertex in $G-\{u,u_1,x\}$ is adjacent to $u$ (and $u$ and $t$ are not adjacnet), we find that $u$ is contained in a triangle of $G-(\{u_1,x,t\}\cup T)$, where $T$ is the set of all $t$-bad vertices in $G-\{u,u_1,x\}$. 
We define $D=\{u_1,x,t\}\cup T$. Lemma \ref{lepljenje} proves that $G-D$ is a WNT. Since $D\subseteq N[t]$, we find that $G$ is reducible. 
This completes the proof of Theorem \ref{reducibilnost}.

 \bigskip

\noindent {\bf Acknowledgement:} 
The author greatly appreciates discussions with Uro\v s Milutinovi${\rm \acute c}$ while working on this paper. 
The author is supported by research grant  P1-0297 of Ministry of Education of Slovenia.


\begin{thebibliography}{99}

\bibitem{campos} C.N.~Campos, Y.~Wakabayashi, On dominating sets of maximal outerplanar graphs, Discrete Appl. Math. 161 (2013), no. 3, 330--335.
\bibitem{diestel} R.~Diestel, Graph theory. Fourth edition. Graduate Texts in Mathematics, 173. Springer, Heidelberg, 2010.
\bibitem{henning2} M.~Dorfling, W.~Goddard, M.A.~ Henning, Domination in planar graphs with small diameter II, Ars Combin. 78 (2006), 237--255.
\bibitem{furuya} M.~Furuya, N.~Matsumoto, A note on the domination number of triangulations, J. Graph Theory 79 (2015), no. 2, 83--85.
\bibitem{henning1} W.~Goddard, M.A,~Henning, Domination in planar graphs with small diameter, J. Graph Theory 40 (2002), no. 1, 1--25.
\bibitem{honjo} T.~Honjo, K.~Kawarabayashi, A.~Nakamoto, Dominating sets in triangulations on surfaces, J. Graph Theory 63 (2010), 17--30.
\bibitem{erika} E.~King, M.~Pelsmajer, Dominating sets in plane triangulations, Discrete Math., 310 (2010),  2221--2230
\bibitem{li} Z.~Li, E.~Zhu, Z.~Shao, J.~Xu, Jin On dominating sets of maximal outerplanar and planar graphs, Discrete Appl. Math. 198 (2016), 164--169.
\bibitem{liu} H.~Liu and M.J.~Pelsmajer, Dominating sets in triangulations on surfaces, Ars Math Contemp. 4 (2011), 177--204.
\bibitem{mac} G.~MacGillivray, K.~Seyffarth, Domination numbers of planar graphs, J. Graph Theory 22 (1996), no. 3, 213--229.
\bibitem{mathesontarjan}  L.R.~Matheson, R.E.~Tarjan, Dominating sets in planar graphs, European J. Combin. 17 (1996), no. 6, 565--568.
\bibitem{mt} B.~Mohar, C.~Thomassen, Graphs on surfaces. Johns Hopkins Studies in the Mathematical Sciences. Johns Hopkins University Press, Baltimore, MD, 2001.
\bibitem{plummer} M.D.~Plummer, X.~Zha, On certain spanning subgraphs of embeddings with applications to domination, Discrete Math., 309 (2009), 4784--4792.
\bibitem{plummer1} M.D.~Plummer, D.~Ye, X.~Zha, Dominating plane triangulations. Discrete Appl. Math., 211 (2016), 175--182 
\bibitem{tokunaga} S.~Tokunaga, Dominating sets of maximal outerplanar graphs, Discrete Appl. Math. 161 (2013), 3097--3099.
\end{thebibliography}
\end{document}